\documentclass[preprint,review,12pt]{elsarticle}
\usepackage[english]{babel}
\usepackage{amsfonts,amssymb,amsmath}
\usepackage{graphicx}

\newtheorem{algorithm}{Algorithm}

\begin{document}

\begin{frontmatter}
\title{Application of the parallel BDDC preconditioner \\
to the Stokes flow}

\author[mu]{Jakub \v{S}\'{\i}stek\corref{cor1}}
\ead{sistek@math.cas.cz}
\author[ucd,ut]{Bed\v{r}ich Soused\'{\i}k}
\ead{bedrich.sousedik@ucdenver.edu}
\author[fs]{Pavel Burda}
\ead{pavel.burda@fs.cvut.cz}
\author[ucd]{Jan Mandel}
\ead{jan.mandel@ucdenver.edu}
\author[fsv]{Jaroslav~Novotn\'{y}}
\ead{novotny@mat.fsv.cvut.cz}

\address[mu]{Institute of Mathematics, Academy of Sciences of the Czech Republic,  \\ 
\v Zitn\' a 25, CZ - 115 67 Prague 1, Czech Republic.}
\address[ucd]{Department of Mathematical and Statistical Sciences,
University of Colorado Denver, \\ 
Denver, CO 80217-3364, USA.}
\address[ut]{Institute~of~Thermomechanics, Academy~of~Sciences of the~Czech~Republic, \\
Dolej\v skova 1402/5, CZ - 182 00  Prague 8, Czech Republic.}
\address[fs]{Department of Mathematics, Faculty of Mechanical Engineering, \\
Czech Technical University in Prague,  \\ 
Karlovo n\'{a}m\v{e}st\'{\i} 13, CZ - 121 35 Prague 2, Czech Republic.}
\address[fsv]{Department of Mathematics,
Faculty of Civil Engineering, \\ Czech Technical University in Prague, 
\\ Th\'{a}kurova 7, CZ - 166 29 Prague 6 , Czech Republic.}
\cortext[cor1]{Corresponding author; 
tel.: +420 222 090 710; fax: +420 222 211 638.}

\begin{abstract}
A parallel implementation of the Balancing Domain Decomposition by Constraints (BDDC) method is described.
It is based on formulation of BDDC with global matrices without explicit coarse problem.
The implementation is based on the MUMPS parallel solver for computing the approximate inverse used for preconditioning.
It is successfully applied to several problems of Stokes flow discretized by Taylor-Hood finite elements 
and BDDC is shown to be a promising method also for this class of problems.
\end{abstract}


\end{frontmatter}

\section{Introduction}

In many areas of engineering, 
numerical solution of problems by the finite element method (FEM) 
leads to solution of systems of linear algebraic equations with sparse and often ill-conditioned matrices.
For very large problems, 
the usual method of choice for their solution is one of the iterative methods based on Krylov subspaces.
However, without preconditioning, the convergence rate deteriorates with growing condition number of the problem.
The need of first-rate preconditioners tailored to the solved problem,
which can be implemented in parallel,
gave rise to the field of domain decomposition methods (e.g. \cite{Toselli-2005-DDM}).

The Balancing Domain Decomposition based on Constraints (BDDC) is one of the most advanced preconditioners of this class.
It was introduced by Dohrmann \cite{Dohrmann-2003-PSC} in 2003 and the theory was developed by Mandel and Dohrmann in \cite{Mandel-2003-CBD}.
In an important contribution to the theory of the preconditioner \cite{Mandel-2005-ATP}, Mandel, Dohrmann, and Tezaur proved 
close connections with the earlier FETI-DP method by Farhat et al.~\cite{Farhat-2001-FDP}, 
another popular domain decomposition technique.
The preconditioner was reformulated without explicit coarse problem as is used in this paper by Li and Widlund in \cite{Li-2006-FBB}.
The underlying theory of the BDDC method covers problems with symmetric positive definite matrix.
An important application that leads to such kind of systems is structural analysis by linear elasticity theory.

The solution of the incompressible Stokes problem by a mixed finite element method leads to a saddle point system with symmetric indefinite matrix.
Thus, the standard theory of BDDC does not cover this important class of problems.
In the first attempt to apply BDDC to the incompressible Stokes problem proposed by Li and Widlund~\cite{Li-2006-BAI}, 
the optimal preconditioning properties of BDDC were recovered.
The approach is based on the notion of \emph{benign subspaces},
which is restricted to using discontinuous pressure approximation,
and the authors present results for piecewise constant functions. 
Moreover, the approach in \cite{Li-2006-BAI} requires quite nonstandard constraints between subdomains,
thus making the implementation more problem specific and difficult.

In this paper, we follow a different approach.
We have implemented a~parallel version of the BDDC method and verified its performance 
on a number of problems arising from linear elasticity (e.g. \cite{Sistek-2010-BFS}).
Here, we investigate the applicability of the method 
and its implementation to the Stokes flow with only minor changes to the source code 
of the implementation for elasticity problems.
Although such application is beyond the standard theory of the BDDC method,
contributive results are obtained.

It has been known for a long time, 
that the conjugate gradient method is able to reach solution also for many indefinite cases (e.g. \cite{Paige-1975-SSI}),
although it may fail in general.
Our effort is also supported by recent trends of numerical linear algebra to investigate and often prefer the use of 
preconditioned CG method (PCG) with a suitable indefinite preconditioner 
over more robust but also more expensive iterative methods for solving indefinite systems 
such as MINRES, BiCG or GMRES \cite{Rozloznik-2002-KSM}.
Another reason for which we do not switch to the MINRES method \cite{Paige-1975-SSI}, 
which is suitable for indefinite problems (e.g. \cite{Elman-2005-FEF}), 
is the fact that it requires a positive definite preconditioner,
while BDDC at the presented setting provides an indefinite preconditioner for the saddle-point problem. 

Several results for the Stokes flow in three dimensions are presented.
All these problems are obtained using mixed discretization by Taylor-Hood finite elements or their serendipity version.
These elements use piecewise (tri)linear pressure approximation,
which does not allow the approach via benign spaces of \cite{Li-2006-BAI},
but are very popular in the computational fluid dynamics community.

\section{Stokes problem and approximation by mixed FEM}
\label{sec:model}

Let $\Omega $ be an open bounded domain in $\mathbb{R}^3$ filled with an incompressible viscous fluid, 
and let $\partial\Omega $ be its boundary. 
Isothermal low speed flow of such fluid is modelled by the following Stokes system of partial differential equations
\begin{eqnarray}
\label{ss}
   - \nu\Delta {\bf u} + \nabla p & = & {\bf f}\ \ \ {\rm in}\ \Omega , \\
\label{rk}   
   - \nabla \cdot {\bf u} & = & 0\ \ \ {\rm in}\ \Omega , \\
\label{bcd}
	{\bf u} & = & {\bf g}\ \ \textrm{on}\ \partial\Omega _g ,\\
\label{bcn}
  - \nu (\nabla {\bf u}) {\bf n} + p{\bf n} & = & {\bf 0}\ \ \textrm{on}\ \partial\Omega _h ,
\end{eqnarray}
where ${\bf u}$ denotes the vector of flow velocity, 
$p$ denotes the pressure divided by the (constant) density, 
$\nu$ denotes the kinematic viscosity of the fluid supposed to be constant,
${\bf f}$ denotes the density of volume forces per mass unit,
$\partial\Omega _g$ and $\partial\Omega _h$ are two subsets of $\partial\Omega $ satisfying 
${\overline{\partial\Omega}} = {\overline {\partial\Omega _g}} \cup {\overline {\partial\Omega _h}},\ \mu_{\mathbb{R}^2} (\partial\Omega_g \cap \partial\Omega _h) = 0$,
${\bf n}$ denotes an outer normal vector to the boundary $\partial\Omega $ with unit length,
and ${\bf g}$ is a given function satisfying $\int_{\partial\Omega } {\bf g}\cdot {\bf n}\ {\rm ds} = 0$ in the case of $\partial\Omega = \partial\Omega _g$ .

We derive the weak formulation of the Stokes equations (\ref{ss})-(\ref{bcn}) in the manner of mixed methods (cf. \cite{Girault-1986-FEM}). 
Let us consider the vector function space 
$V= \Big\{ {\bf v}\in [H^1(\Omega )]^3; {\bf v}|_{\partial\Omega _g} = {\bf 0} \Big\}$ 
and the set 
$V_g = \Big\{ {\bf v}\in [H^1(\Omega )]^3;{\bf v}|_{\partial\Omega _g} = {\bf g} \Big\}$,
where $H^1(\Omega )$ is the usual Sobolev space, 
and the restriction ${\bf v}|_{\partial\Omega _g}$ is understood in the sense of traces.

We now introduce a triangulation of the domain $\Omega $
into Taylor-Hood finite elements $P_{2}P_{1}$ and/or $Q_{2}Q_{1}$ (or their serendipity version $Q_{2S}Q_{1}$), 
which satisfy the Babu\v ska-Brezzi stability condition (cf. \cite{Brezzi-1991-MHF}).
Their application leads to the finite dimensional subsets 
$V_{gh} \subset V_g$, $V_h \subset V$, 
which contain continuous piecewise quadratic functions,
and $Q_h \subset L_2(\Omega )$ with continuous piecewise linear functions.

We can now introduce the \emph{discrete Stokes problem}:

\noindent Find ${\bf u}_h \in V_{gh}$ and $p_h \in Q_h$ satisfying
\begin{eqnarray}
\nonumber
\label{ssh}  
  \nu\int_{\Omega } \nabla {\bf u}_h :\nabla {\bf v}_h {\rm d}\Omega - \int_{\Omega } p_h \nabla \cdot {\bf v}_h {\rm d}\Omega &=& \int_{\Omega } {\bf f} \cdot {\bf v}_h {\rm d}\Omega ,\ \forall {\bf v}_h \in V_h,\\
\label{rksemih}  
  - \int_{\Omega } \psi _h\nabla \cdot {\bf u}_h {\rm d}\Omega &=& 0,\ \forall \psi _h \in Q_h,\\
\label{bcwsemih}  
  {\bf u}_h - {\bf u}_{gh} &\in & V_h.
\end{eqnarray}
Here ${\bf u}_{gh} \in V_{gh}$ represents the Dirichlet boundary condition ${\bf g}$ in (\ref{bcd}).

Expressing the finite element functions as linear combinations of basis functions (see e.g. \cite[Section 5.3]{Elman-2005-FEF} for more details),
the problem finally leads to the saddle point system of algebraic equations
\begin{equation}
\label{eq:saddle-point}
\left[
\begin{array}{cc}
\mathbf{A} & B^T \\
B & 0
\end{array}
\right]
\left[
\begin{array}{c}
\mathbf{\bar{u}} \\
\mathbf{\bar{p}}  
\end{array}
\right]
=
\left[
\begin{array}{c}
\mathbf{\bar{f}} \\
\mathbf{0}  
\end{array}
\right] ,
\end{equation}
where $\mathbf{\bar{u}}$ denotes velocity unknowns, $\mathbf{\bar{p}}$ denotes pressure unknowns, ${\bf A}$ and ${B}$ are called vector--Laplacian matrix and
the divergence matrix, respectively, and ${\bf \bar{f}}$ is the discrete vector of intensity of volume forces per mass unit.

\section{Iterative substructuring}
\label{sec:substructuring}

In this section, we recall ideas of iterative substructuring used in our implementation.
Details can be found in \cite{Toselli-2005-DDM}.

Let $\Omega $ be a~bounded domain in $\mathbb{R}^{2}$ or $\mathbb{R}^{3}$,
let $U$ be a~finite element space of piecewise polynomial functions $v$ continuous on $\Omega$ and $U'$ its dual space.
Let $a(\cdot,\cdot)$ be a~bilinear form on $U \times U$ and $f\in U'$,
and let $\langle \cdot,\cdot \rangle$ denote the duality pairing of $U'$ and $U$.
Consider now an abstract variational problem:
\emph{Find $u\in U$ such that}
\begin{equation}
\label{eq:var}
a(u,v)=\langle f,v\rangle\quad\forall\,v\in U\,.
\end{equation}

Write the matrix problem corresponding to \eqref{eq:var} as
\begin{equation}
\label{eq:A-problem}
Au=f.
\end{equation}

The domain $\Omega$ is decomposed into $N$ nonoverlapping subdomains $\Omega_{i}$, $i=1,...,N$, with characteristic size $H$,
which form a~conforming triangulation of the domain $\Omega $.
Each subdomain is a~union of several finite elements of the underlying mesh with characteristic mesh size $h$,
i.e. nodes of the finite elements between subdomains coincide.
Unknowns common to at least two subdomains are called \emph{boundary unknowns} and the union of all
boundary unknowns is called the \emph{interface} $\Gamma$.

The problem is first reduced to the interface $\Gamma$.
For this purpose, the solution $u$ is split into the 
interior solution $u_{o}$, with zero values at $\Gamma $, and $u_{\Gamma}$, 
where values in subdomain interiors are determined by values at $\Gamma $ (see (\ref{eq:discharm}) below). 
Then problem (\ref{eq:A-problem}) may be rewritten as
\begin{equation}
A(u_{\Gamma}+u_{o})=f. \label{eq:A-problem2}
\end{equation}
Let us formally reorder unknowns of problem (\ref{eq:A-problem2}) into two blocks, with the
first block (subscript $_1$) corresponding to unknowns in subdomain interiors, and the second
block (subscript $_2$) corresponding to unknowns at the interface.
This results in the block form of the system (\ref{eq:A-problem2}) given as
\begin{equation}
\left[
\begin{array}
[c]{cc}
A_{11} & A_{12}\\
A_{21} & A_{22}
\end{array}
\right]  \left[
\begin{array}
[c]{c}
u_{\Gamma1}+u_{o1}\\
u_{\Gamma2}+u_{o2}
\end{array}
\right]  =\left[
\begin{array}
[c]{c}
f_{1}\\
f_{2}
\end{array}
\right]  , \label{eq:A-block}
\end{equation}
with $u_{o2}=0$ by definition. 
Function $u_{\Gamma}$ is called \emph{discrete harmonic},
by which we mean that it is fully determined by values at interface $\Gamma $ and by the algebraic condition 
\begin{equation}
\label{eq:discharm}
A_{11} u_{\Gamma1} + A_{12} u_{\Gamma2} = 0,\ \ \mbox{i.e.} \ \ A_{11} u_{\Gamma1} = - A_{12} u_{\Gamma2}.
\end{equation}
By this splitting, we derive that
(\ref{eq:A-block}) is equivalent to
\begin{equation}
A_{11}u_{o1}=f_{1}, \label{eq:A-interior}
\end{equation}
\begin{equation}
\left[
\begin{array}
[c]{cc}
A_{11} & A_{12}\\
A_{21} & A_{22}
\end{array}
\right]  \left[
\begin{array}
[c]{c}
u_{\Gamma1}\\
u_{\Gamma2}
\end{array}
\right]  =\left[
\begin{array}
[c]{c}
0\\
f_{2}-A_{21}u_{o1}
\end{array}
\right]  , \label{eq:A-block-interface2}
\end{equation}
and the solution is obtained as $u=u_{\Gamma}+u_{o}$. 
Problem
(\ref{eq:A-block-interface2}) 
can be further split into two problems
\begin{equation}
A_{11}u_{\Gamma1}=-A_{12} u_{\Gamma2}, \label{eq:S-interior}
\end{equation}
\begin{equation}
Su_{\Gamma2}=g_{2}, \label{eq:S-problem}
\end{equation}
where $S$ is the \emph{Schur complement} with respect to unknowns at interface $\Gamma$ defined as
$S=A_{22}-A_{21}A_{11}^{-1}A_{12}$, and $g_{2}$ is the \emph{condensed right
hand side} $g_{2}=f_{2}-A_{21}u_{o1} = f_{2}-A_{21}A_{11}^{-1}f_{1}$. 
Since $A_{11}$ has a block diagonal structure, the solution to
(\ref{eq:A-interior}) may be found in parallel and similarly the solution to
(\ref{eq:S-interior}).
We are ready to recall the algorithm of substructuring.

\begin{algorithm}[Iterative substructuring]
\label{alg:substructuring}
Problem (\ref{eq:A-problem}) is solved in the following steps:
\begin{enumerate}
\item factorize block diagonal matrix $A_{11}$ in (\ref{eq:A-interior}) and store factors,
\item solve (\ref{eq:A-interior}) by back-substitution to find $u_{o1}$,
\item construct $g_{2}$ as $g_{2}=f_{2}-A_{21}u_{o1}$,
\item solve problem (\ref{eq:S-problem}) by a Krylov subspace method. 
In each iteration, multiplication of a given vector $p_2$ by $S$ is realized as 
\begin{itemize}
\item find $p_1$ by solution of $A_{11}p_{1}= -A_{12}p_{2}$,
\item get $Sp_2$ as $Sp_2 = A_{21}p_{1} + A_{22}p_{2}$.
\end{itemize}
\item Find $u_{\Gamma 1}$ by (\ref{eq:S-interior}),
\item get solution $u$ as $u = \left[
\begin{array}
[c]{c}
u_{\Gamma1}\\
u_{\Gamma2}
\end{array}
\right]
+ \left[
\begin{array}
[c]{c}
u_{o1}\\
0
\end{array}
\right]
$.
\end{enumerate}
\end{algorithm}
Note, that the Schur complement $S$ is never formed explicitly and its action is realized by three sparse matrix multiplications 
and one back-substitution.
The main reason for using Algorithm \ref{alg:substructuring} is usually much faster convergence of the iterative method for problem (\ref{eq:S-problem}) 
compared to problem (\ref{eq:A-problem}) (see e.g. \cite{Toselli-2005-DDM}).

\section{BDDC preconditioner}
\label{sec:bddc}

The BDDC method provides a preconditioner for problem (\ref{eq:S-problem}).
Let $W_{i}$ be the space of finite element functions on subdomain $\Omega _{i}$, 
coefficients of which satisfy the algebraic discrete harmonic condition 
(\ref{eq:discharm}) locally on the subdomain,
and put $W=W_{1}\times \cdots \times W_{N}$.
It is the space where subdomains are completely disconnected at the interface $\Gamma $,
 and functions on them are independent of each other.
Let us further define $U_{\Gamma } \subset U$ as the subset of finite element functions on $\Omega $ with coefficients satisfying 
the discrete harmonic condition (\ref{eq:discharm}).
Clearly, $U_{\Gamma } \subset W$, and the solution to (\ref{eq:A-block-interface2}) $u_{\Gamma} \in U_{\Gamma }$.

The main idea of the BDDC preconditioner
in the abstract form \cite{Mandel-2007-BFM} is to construct an auxiliary finite
dimensional space $\widetilde{W}$ such that
$U_{\Gamma } \subset\widetilde{W}\subset W$,
and extend the bilinear form $a\left(  \cdot,\cdot\right)  $ to a~form
$\widetilde{a}\left(  \cdot,\cdot\right) $ defined on $\widetilde{W}
\times\widetilde{W}$, such that solving the variational problem
\eqref{eq:var} with $\widetilde{a}\left(  \cdot,\cdot\right)  $ in place of
$a\left(  \cdot,\cdot\right)  $ is cheaper and can be split into
independent computations performed in parallel. Then the solution restricted to
$U_{\Gamma } $ is used for the preconditioning of \eqref{eq:S-problem}.
Space $\widetilde{W}$ contains functions generally discontinuous at interface $\Gamma $ 
except a small set of \emph{coarse degrees of freedom} at which continuity is preserved.
Coarse degrees of freedom are typically values at selected nodes called \emph{corners}.
In addition, 
continuity of generalized degrees of freedom, such as \emph{averages} over subdomain \emph{edges} 
and/or \emph{faces}, might be enforced.

\begin{figure}[htbp]
   \begin{center}
   \includegraphics[width=60mm]{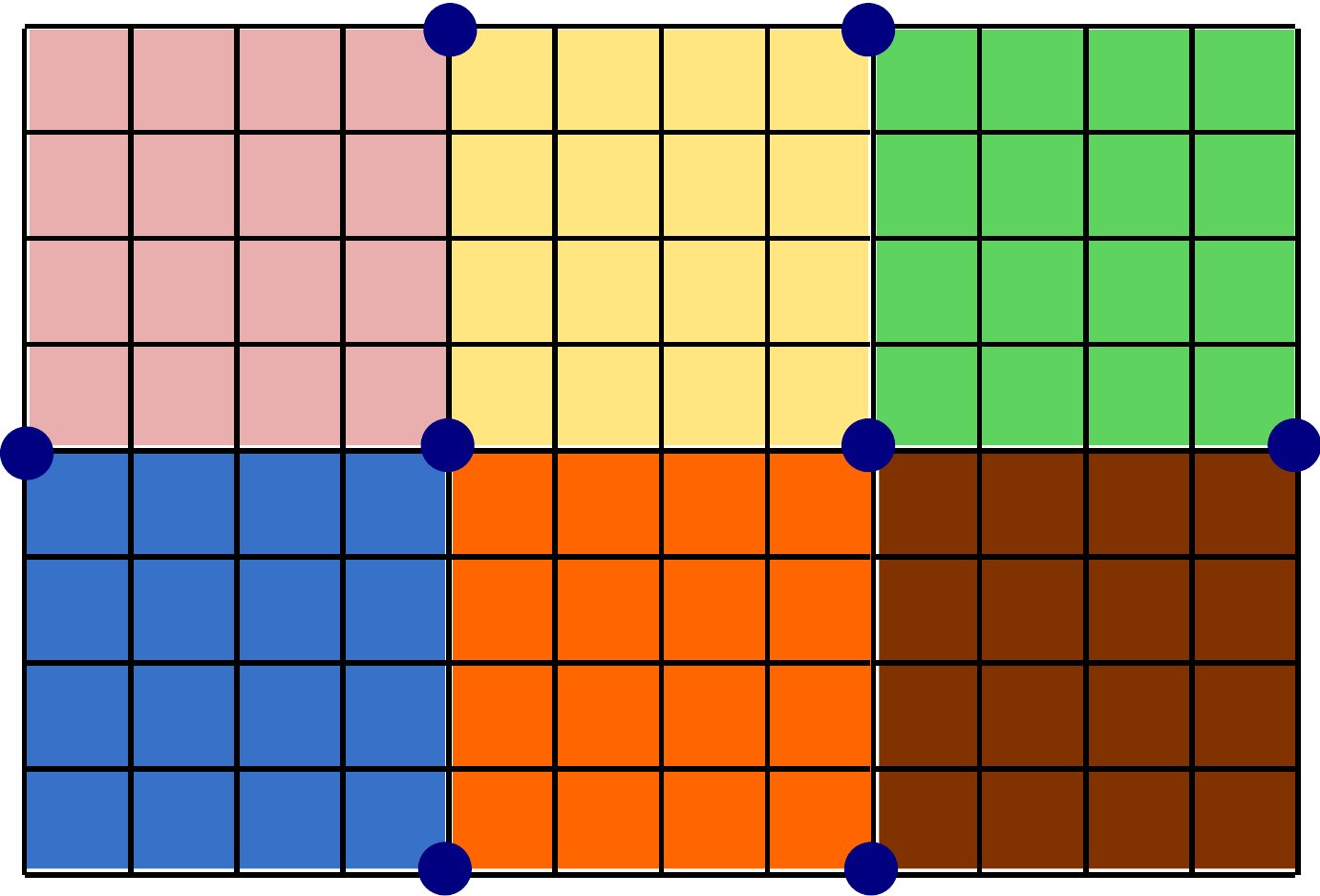} 
   \ \
   \includegraphics[width=60mm]{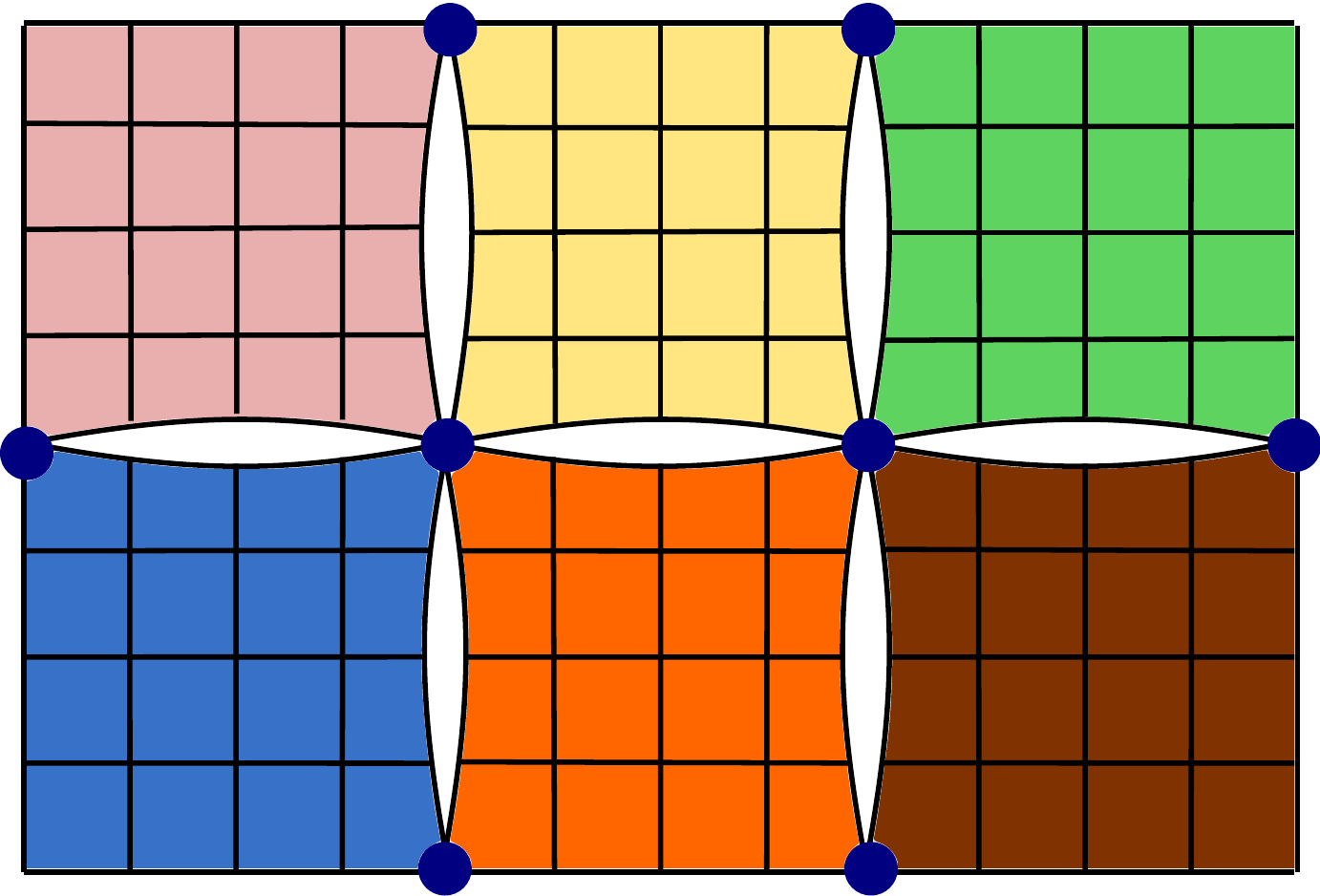}
   \caption{\label{fig:virtual_mesh} 
      Example of an actual computational mesh with six subdomains (left) and corresponding \emph{virtual mesh} used for assembly of matrix $\widetilde{A}$ (right).}
   \end{center}
\end{figure}

In computation, the corresponding matrix denoted $\widetilde{A}$ is used.
It is larger than the original matrix of the problem $A$,
but it possesses a simpler structure suitable for direct solution methods.
This is the reason why it can be used as a preconditioner.
In the presented algorithm, 
matrix $\widetilde{A}$ is constructed using the standard FEM assembly procedure on a~virtual mesh 
which is disconnected at interface outside corners (Figure \ref{fig:virtual_mesh}).

The projection $E:\widetilde{W}\rightarrow U_{\Gamma } $ is realized as a~weighted average of values from different
subdomains at unknowns on the interface $\Gamma$, thus resulting in functions
continuous across the interface.
The weights at a degree of freedom are chosen as the inverse to the number of subdomains, 
in which the degree of freedom is contained, as is done e.g. in \cite{Li-2006-BAI}. 
This approach is used for both velocity and pressure unknowns.

Let $r \in U_{\Gamma }'$ be the residual in an iteration of an iterative method.
The BDDC preconditioner $M_{BDDC}: U_{\Gamma }'\rightarrow U_{\Gamma }$ in the abstract form (see \cite{Mandel-2007-BFM}) produces the preconditioned residual $v \in U_{\Gamma }$ as
\begin{equation*}
M_{BDDC}: r\rightarrow v=Ew,
\end{equation*}
where $w\in \widetilde{W}$ is obtained as the solution to problem
\begin{equation}
w\in\widetilde{W}:\widetilde{a}\left(  w,z\right)  =\left(  r,Ez\right) \quad\forall z\in\widetilde{W}, \label{eq:bddc}
\end{equation}
or in terms of matrices as
\begin{equation}
\label{eq:bddcmatrix}
v=E\widetilde{A}^{-1}E^{T}r.
\end{equation}
Here, the action of $\widetilde{A}^{-1}$ is performed as a back-substitution by a direct solver.
After $v$ is found, we are typically interested only in its values at interface nodes, 
since multiplication of this vector by $S$ follows and interior values are resolved 
from discrete harmonic constraint (\ref{eq:discharm}) as described in Algorithm \ref{alg:substructuring}.

For the Stokes problem, we adopt the following slightly unusual notation in (\ref{eq:var})
\begin{eqnarray}
a(u,v)&=&\nu\int_{\Omega } \nabla {\bf u}_h : \nabla {\bf v}_h {\rm d}\Omega - \int_{\Omega } p_h \nabla \cdot {\bf v}_h {\rm d}\Omega - \int_{\Omega } \psi _h\nabla \cdot {\bf u}_h {\rm d}\Omega , \label{eq:astokes} \\
\langle f,v\rangle&=&\int_{\Omega} {\bf f} \cdot {\bf v}_h {\rm d}\Omega. \label{eq:rhsstokes}
\end{eqnarray}
The bilinear form $a(u,v)$ is symmetric but indefinite \cite{Brezzi-1991-MHF,Elman-2005-FEF}.
In (\ref{eq:A-problem}), this corresponds to putting 
\begin{equation*}
A = \left[
\begin{array}{cc}
\mathbf{A} & B^T \\
B & 0
\end{array}
\right],
u = 
\left[
\begin{array}{c}
\mathbf{\bar{u}} \\
\mathbf{\bar{p}}  
\end{array}
\right],
f =
\left[
\begin{array}{c}
\mathbf{\bar{f}} \\
\mathbf{0}  
\end{array}
\right] .
\end{equation*}

Should we distinguish between blocks corresponding to nodes in interiors of subdomains and at the interface $\Gamma$ as in (\ref{eq:A-block}),
saddle point problem (\ref{eq:saddle-point}) would look as
\begin{equation}
\left[
\begin{array}{cccc}
\mathbf{A}_{11} & \mathbf{A}_{12} & B^T_{11} & B^T_{21} \\
\mathbf{A}_{21} & \mathbf{A}_{22} & B^T_{12} & B^T_{22} \\
        B_{11}  &         B_{12}  & 0        & 0        \\
        B_{21}  &         B_{22}  & 0        & 0        
\end{array}
\right]
\left[
\begin{array}{c}
\mathbf{\bar{u}_1} \\
\mathbf{\bar{u}_2} \\
\mathbf{\bar{p}_1} \\
\mathbf{\bar{p}_2}  
\end{array}
\right]
=
\left[
\begin{array}{c}
\mathbf{\bar{f}_1} \\
\mathbf{\bar{f}_2} \\
\mathbf{0} \\ 
\mathbf{0}
\end{array}
\right] ,
\end{equation}
and the Schur complement matrix and the condensed right hand side in (\ref{eq:S-problem}) as 
\begin{eqnarray*}
S &=& 
\left[
\begin{array}{cc}
\mathbf{A}_{22} & B^T_{22} \\
        B_{22}  & 0        
\end{array}
\right] 
- 
\left[
\begin{array}{cc}
\mathbf{A}_{21} & B^T_{12}  \\
        B_{21}  & 0         
\end{array}
\right]
\left[
\begin{array}{cc}
\mathbf{A}_{11} & B^T_{11}  \\
        B_{11}  & 0        
\end{array}
\right]^{-1}
\left[
\begin{array}{cc}
 \mathbf{A}_{12} & B^T_{21} \\
         B_{12}  & 0 
\end{array}
\right], \\
g_2 &=& 
\left[
\begin{array}{c}
\mathbf{\bar{f}_2} \\
\mathbf{0}
\end{array}
\right]
- 
\left[
\begin{array}{cc}
\mathbf{A}_{21} & B^T_{12}  \\
        B_{21}  & 0         
\end{array}
\right]
\left[
\begin{array}{cc}
\mathbf{A}_{11} & B^T_{11}  \\
        B_{11}  & 0        
\end{array}
\right]^{-1}
\left[
\begin{array}{c}
\mathbf{\bar{f}_1} \\
\mathbf{0}
\end{array}
\right] .
\end{eqnarray*}

\section{Implementation and numerical results}
\label{sec:numericalresults}

Our parallel implementation of the BDDC preconditioner has been extensively tested on problems with symmetric positive definite matrices
arising from linear elasticity (e.g. \cite{Sistek-2010-BFS}).
The current version is based on the multifrontal massively parallel sparse direct solver MUMPS \cite{Amestoy-2000-MPD} version 4.9.2,
which is used for the factorization of the matrices $\widetilde{A}$ in (\ref{eq:bddcmatrix}) and $A_{11}$ in (\ref{eq:S-interior}).
These matrices are put into MUMPS in the distributed format with one subdomain corresponding to one processor.
While Cholesky factorization is used for problems with symmetric positive definite system matrices, 
for the Stokes problem, MUMPS is simply switched to the $LDL^T$ factorization of general symmetric matrices.
If additional constraints on averages over edges or faces are prescribed, 
the generalized change of variables is used in combination with nullspace projection.
This approach, which provides a generalization of the change of basis from \cite{Li-2006-FBB}, is described in detail in \cite{Mandel-2009-ABT}.
Iterations are performed by a parallel PCG solver.

The applicability of the preconditioner to the steady problem of Stokes flow has been tested,
and results are presented in this section.
The system matrix of the Stokes problem is symmetric, but indefinite.
For this reason, the standard theory of BDDC does not cover this case.
A way to assure positive definiteness of the preconditioned operator based on BDDC was presented by Li and Widlund \cite{Li-2006-BAI}.
However, that approach is limited to discontinuous pressure approximation, 
and thus it can be used for neither $Q_{2}Q_{1}$ Taylor--Hood finite elements (e.g. \cite{Brezzi-1991-MHF}) used in our Matlab computations 
nor their serendipity version $Q_{2S}Q_{1}$ used in our parallel computations.

\subsection{Problem (1)}

\begin{figure}[htbp]
   \begin{center}
   \includegraphics[width=50mm]{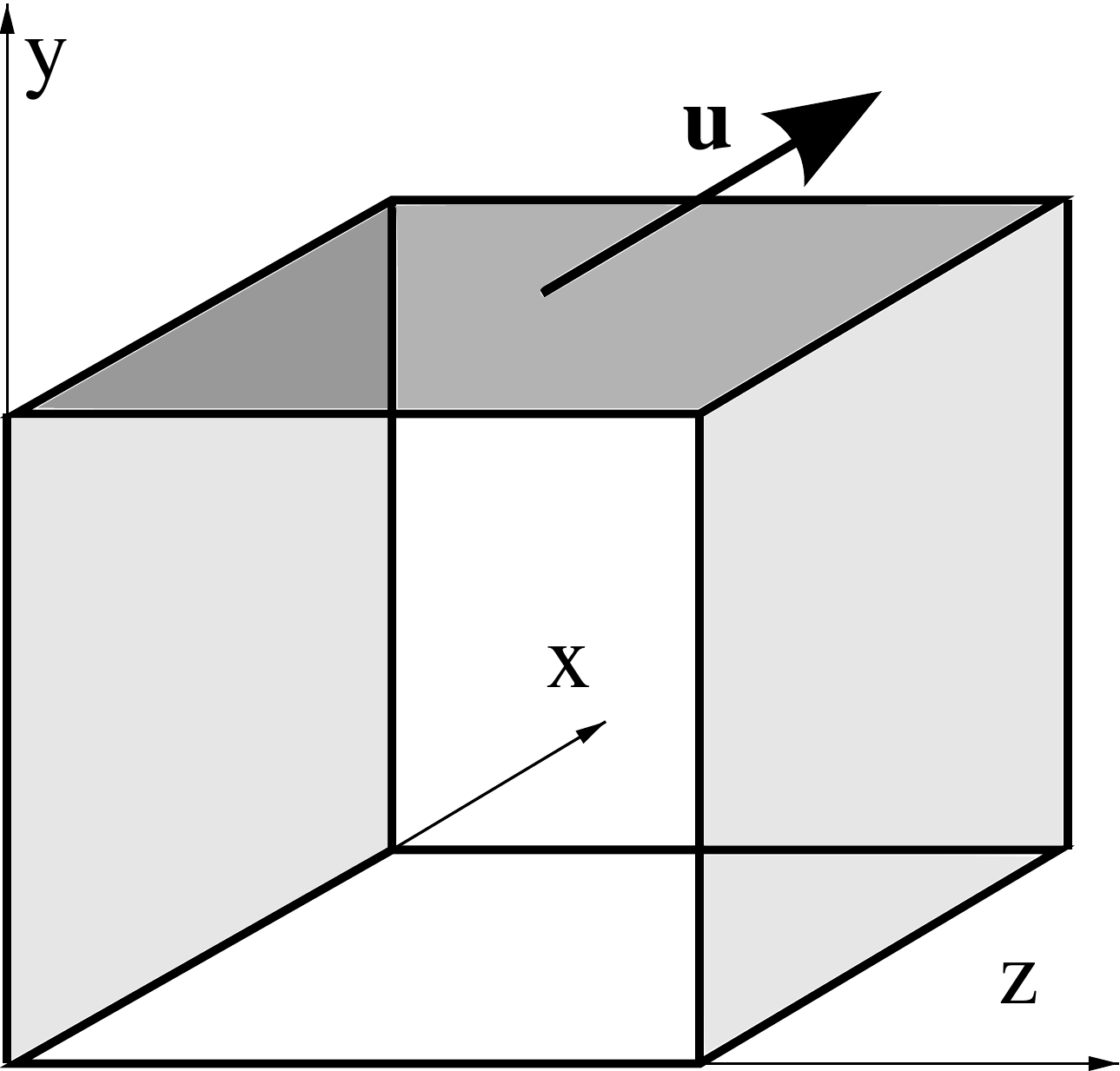}
   \hskip 20mm
   \includegraphics[width=50mm]{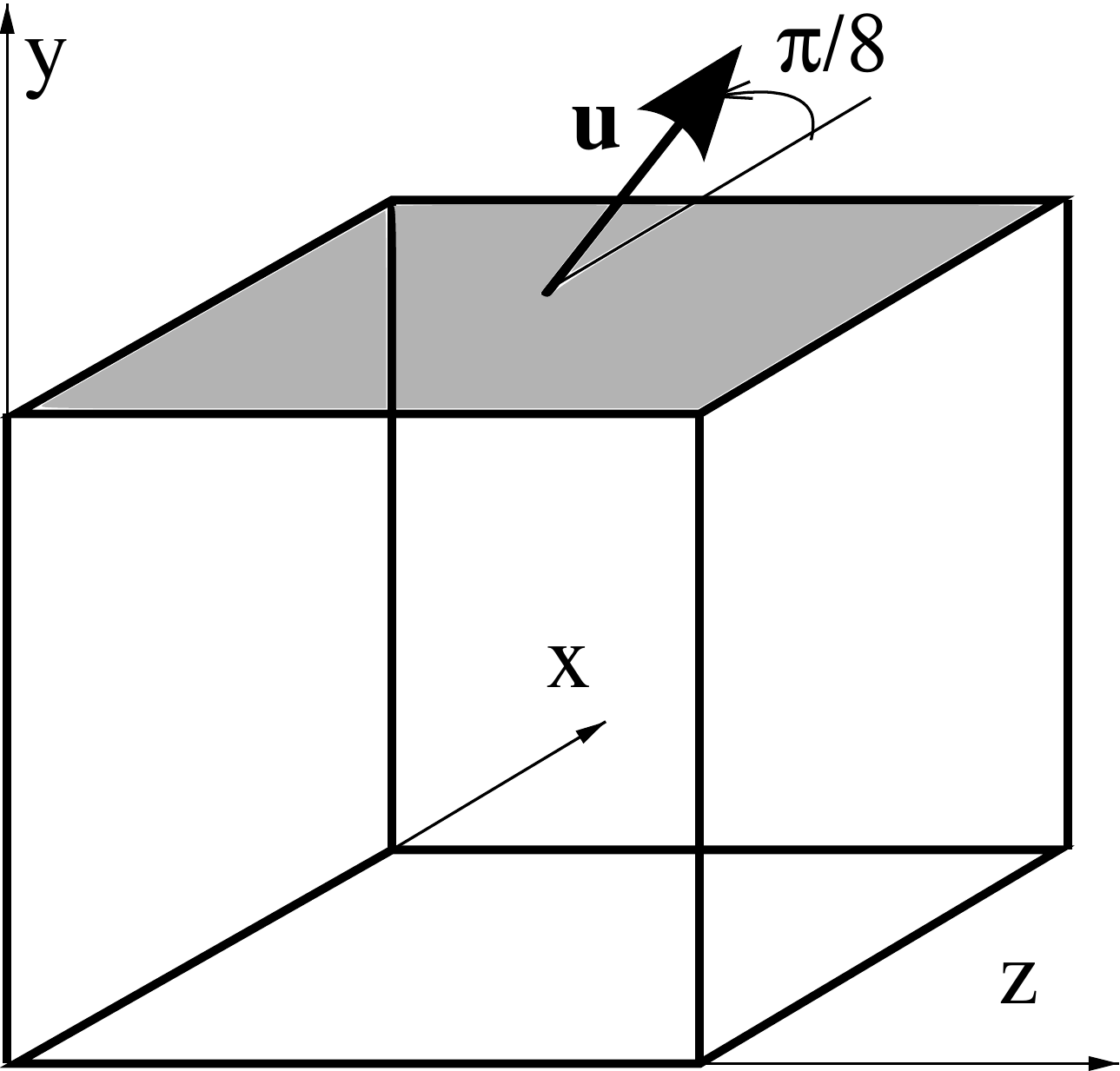}
   \caption{\label{fig:kavita3d} 
      3D lid driven cavity problems: 
      section of infinite cavity (Problem 1) (left) -- unit velocity aligned with $x$-axis prescribed on the upper (dark grey) face, zero component of velocity in $z$ direction on two (light grey) faces, homogeneous Dirichlet boundary conditions for velocity on remaining (white) faces;
      cubic cavity (Problem 2) (right) -- unit velocity parallel to the $xz$-plane rotated along $y$-axis prescribed on the upper (dark grey) face, homogeneous Dirichlet boundary condition for velocity on remaining (white) faces.
   }
   \end{center}
\end{figure}

The method is first tested on the problem of the lid driven cavity.
This popular 2D benchmark problem is used in the 3D setting as a section of an~infinite cavity (Figure \ref{fig:kavita3d} left) 
used e.g. in \cite{Gartling-2006-QFE}:
The domain is a~unit cube with unit velocity in the direction of the $x$-axis on the upper face (called \emph{lid}),
zero normal component of velocity ($u_z = 0$) prescribed on faces parallel to $xy$-plane,
and homogeneous Dirichlet boundary conditions for velocity on remaining faces.
In numerical solution, all nodes with $y = 1$ are included into the lid, 
which means that the setting corresponds to the so called \emph{leaky} cavity \cite{Elman-2005-FEF}.
We fix pressure at the node in the centre of the domain to make its solution unique.
The entire motion inside cavity is driven by viscosity of the fluid which is chosen as 0.01.

The problem is uniformly discretized using $Q_{2S}Q_{1}$ serendipity finite elements (velocity unknowns at vertices and edge-centres, pressure unknowns at vertices).

Tables \ref{tab:scaling_cavityHh4}--\ref{tab:scaling_cavityHh16} contain results for variable number of subdomains (columns) and variable $H/h$ ratio,
where $H$ stands for the characteristic size of a subdomain and $h$ denotes the characteristic size of an element.
Each column contains results for constraints at corners only (`c') and with additional constraints on averages over edges (`c+e').
Results are summarised with respect to the number of PCG iterations and computational times of individual dominant operations in the preconditioner 
(the total time includes also time for parallel factorization of the block of interior unknowns).
Computations were performed on SGI Altix~4700 computer in CTU Supercomputing Centre in Prague with 72 1.5~GHz Intel Itanium~2 processors.
The stopping criterion for PCG was defined by relative residual as $\|r\|_{2}/\|g\|_{2} < 10^{-6}$.

Resulting times are for selected settings compared to those by direct application of the MUMPS solver \cite{Amestoy-2000-MPD} and to results by our in--house serial direct solver.
The latter is based on the unsymmetric frontal method by Hood \cite{Hood-1976-FSP},
which is a generalization of the classical frontal method  developed for symmetric positive definite problems by Irons \cite{Irons-1970-FSS}.
The basic idea of the frontal approach lies in simultaneous assembling and eliminating of rows of the system matrix finding the factors `out-of-core'. 
For suitably numbered elements, this approach usually results in huge reduction of memory requirements, 
while it might have a negative impact on the speed due to the I/O operations inside factorization.
This solver has been successfully used by our group to solve a number of benchmark as well as real-life Stokes and Navier-Stokes problems.

Solution of this cavity problem using 16$\times$16$\times$16 elements 
(corresponding to the last column of Table \ref{tab:scaling_cavityHh4} or the first column of Table \ref{tab:scaling_cavityHh8}) 
takes 30 minutes by our frontal solver on a~single processor.

Table \ref{tab:scaling_cavityHh16} also summarises solution of the problem using  32$\times$32$\times$32 elements divided into 8 subdomains
by the MUMPS solver.
We have not been able to fit this problem into memory using the serial frontal solver.

\begin{figure}[htbp]
   \begin{center}
   \includegraphics[width=60mm]{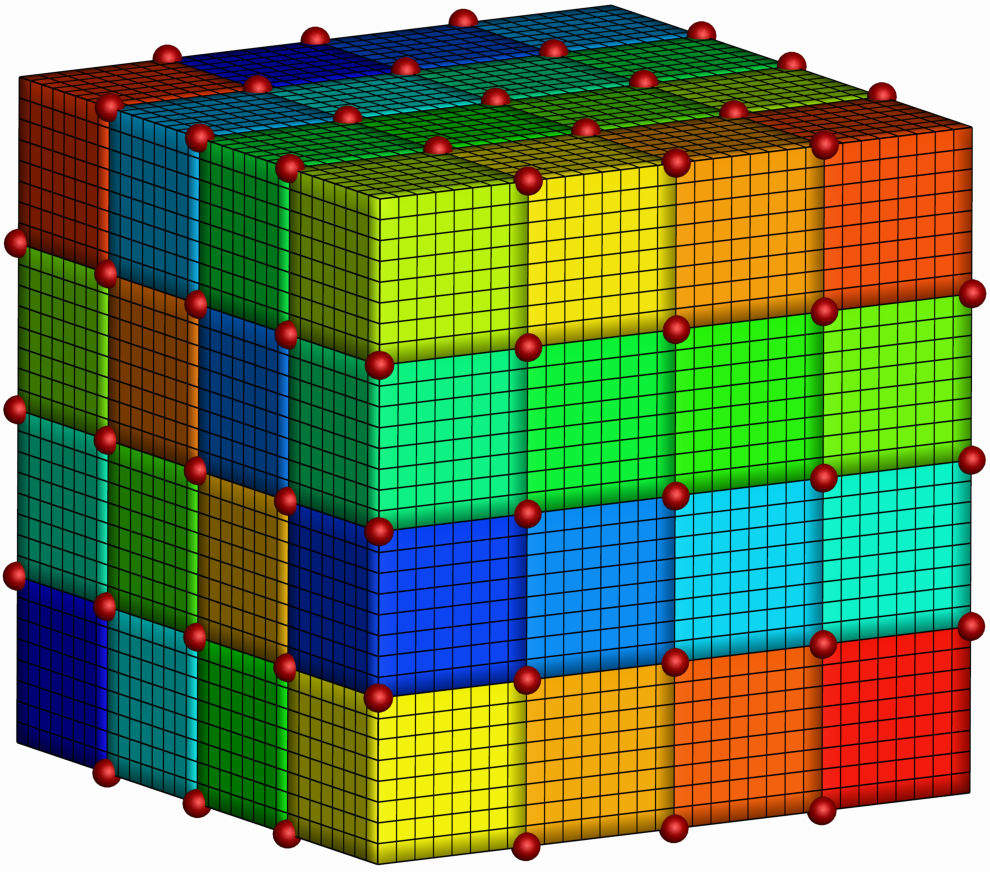}
   \includegraphics[width=70mm]{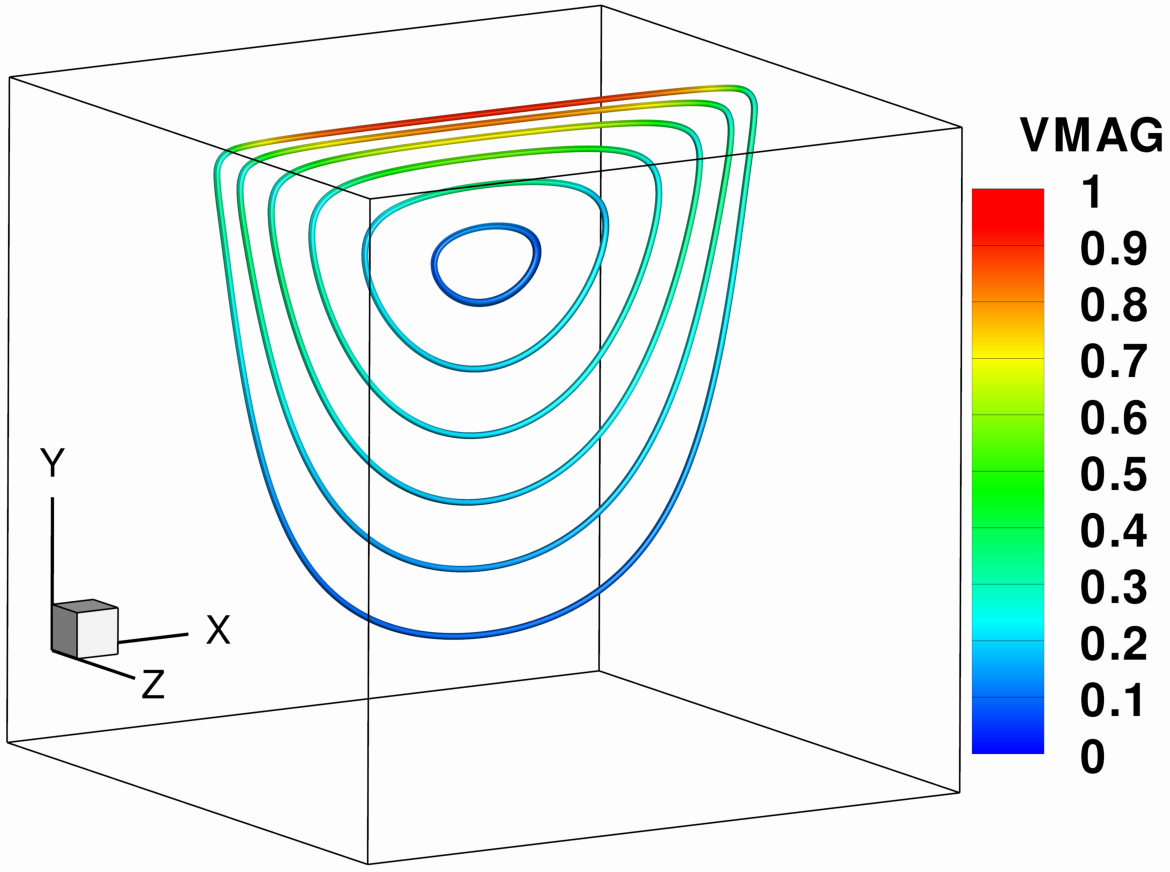}
   \caption{\label{fig:problem_1_4} 3D lid driven cavity problem with 32$\times $32$\times $32 elements, $H/h = 8$: 
                                    division into 64 regular subdomains (left); 
                                    several streamlines at the $z = 0.5$ plane for Problem (1) (right), colours by velocity magnitude.}
   \end{center}
\end{figure}

An example of the problem for 32$\times $32$\times $32 = 32,768 elements and $H/h = 8$ (64 subdomains) is 
presented in Figure \ref{fig:problem_1_4} left. 
In Figure \ref{fig:problem_1_4} right, several streamlines at the $z = 0.5$ plane are presented. 
These are coloured by the velocity magnitude.

\begin{table}[htbp]
\centering
\begin{tabular}{|c|c|c|c|}
\hline
subdomains (processors)      &      8      &      27     &      64     \\ 
unknowns                     &    8,748    &    27,040   &    61,268   \\ \hline
constraints                  & c    / c+e  & c    / c+e  & c    / c+e  \\ \hline
number of PCG iterations     & 18   / 15   & 29   / 17   & 37   / 18   \\ \hline
analysis by MUMPS (sec)      & 0.9  / 0.5  & 1.9  / 2.1  & 5.7  / 6.6  \\
factorization by MUMPS (sec) & 0.2  / 0.2  & 0.3  / 0.3  & 0.5  / 0.6  \\
PCG iterations (sec)         & 2.1  / 1.8  & 12   / 7.1  & 46   / 23   \\ \hline
one PCG iteration (sec)      & 0.12 / 0.12 & 0.42 / 0.42 & 1.24 / 1.26 \\ \hline
total wall time (sec)        & 5.6  / 4.9  & 24   / 22   & 106  / 87   \\ \hline
\end{tabular}
\caption{Scaling of cavity problem (1) for variable number of subdomains and types of constraints, $H/h = 4$.}
\label{tab:scaling_cavityHh4}
\end{table}

\begin{table}[htbp]
\centering
\begin{tabular}{|c|c|c|c|c|}
\hline
subdomains (processors)      &      8      &      27      &      64     \\ 
unknowns                     &   61,268    &    197,500   &     457,380 \\ \hline
constraints                  & c    / c+e  &  c    / c+e  & c    / c+e  \\ \hline
number of PCG iterations     & 19   / 16   &  35   / 18   & 122  / 44   \\ \hline
analysis by MUMPS (sec)      & 11   / 11   &  26   / 27   & 84   / 89   \\
factorization by MUMPS (sec) & 1.3  / 1.4  &  1.9  / 2.1  & 3.0  / 3.0  \\
PCG iterations (sec)         & 9.7  / 8.4  &  70   / 37   & 724  / 281  \\ \hline
one PCG iteration (sec)      & 0.51 / 0.52 &  2.01 / 2.08 & 5.93 / 6.39 \\ \hline
total wall time (sec)        & 42   / 41   &  195  / 182  & 1,352 / 966 \\ \hline
\end{tabular}
\caption{Scaling of cavity problem (1) for variable number of subdomains and types of constraints, $H/h = 8$.}
\label{tab:scaling_cavityHh8}
\end{table}

\begin{table}[htbp]
\centering
\begin{tabular}{|c|c||c|}
\hline
subdomains (processors)      & \multicolumn{2}{|c|}{8}       \\ 
unknowns                     & \multicolumn{2}{|c|}{457,380} \\ \hline
solver                       & BDDC + PCG  &     MUMPS     \\ \hline
constraints                  & c    / c+e  &      n/a      \\ \hline
number of PCG iterations     & 45   / 36   &      n/a      \\ \hline
analysis by MUMPS (sec)      & 168  / 166  &      185      \\
factorization by MUMPS (sec) & 54   / 55   &     1,398     \\
PCG iterations (sec)         & 125  / 110  &      n/a      \\ \hline
one PCG iteration (sec)      & 2.78 / 3.05 &      n/a      \\ \hline
total wall time (sec)        & 634  / 651  &     1,601     \\ \hline
\end{tabular}
\caption{Cavity problem (1): BDDC for variable types of constraints, and direct solution by MUMPS solver, $H/h = 16$.}
\label{tab:scaling_cavityHh16}
\end{table}

Tables~\ref{tab:scaling_cavityHh4}--\ref{tab:scaling_cavityHh16} reveal an unfortunate property of the presented approach, 
that scalability is not achieved with respect to number of resulting iterations. 
Computational times are growing with growing number of subdomains not only due to increasing number of iterations, 
but also due to the dependence on the MUMPS solver, which turns out not to scale well for this problem.

Nevertheless, it is still interesting to compare the results to the frontal solver used to address the Stokes and Navier--Stokes problems by our group before, 
and even to compare the computational time of solution by PCG with BDDC preconditioner and by MUMPS (Table~\ref{tab:scaling_cavityHh16}) on eight processors, 
for which the former is two times faster.
This result supports the initial idea of the implementation -- that using the parallel direct solver for the disconnected problem as a preconditioner for 
an iterative method can be much faster than using the same solver for the original problem directly.

\subsection{Problem (2)}

The next problem is formed by another generalization of the 2D cavity problem into 3D inspired by \cite{Wathen-2003-NPO}.
It has the following setting (Figure \ref{fig:kavita3d} right):
Homogeneous Dirichlet boundary conditions for velocity are considered on all faces except the lid, 
where unit tangential velocity is prescribed.
However, to make also the nature of the flow three-dimensional, the velocity is now not aligned with the $x$-axis but rotated by angle $\pi /8$.
Again, pressure is fixed at the node in the centre of the domain and the viscosity of the fluid is chosen as 0.01.

This problem is uniformly discretized using (full) $Q_{2}Q_{1}$ Taylor--Hood finite elements.
Serendipity $Q_{2S}Q_{1}$ elements cannot be used for this problem since they fail to determine pressure at the eight corners of the domain
if Dirichlet boundary conditions on velocity are prescribed on all three faces at a~corner.

Tables~\ref{tab:iter_cavityHh4}--\ref{tab:iter_cavity222} summarise performance of the BDDC preconditioner in connection with GMRES and BiCGStab methods, respectively.
Since these experiments were run in Matlab by a serial implementation of the method, the comparison is done only for number of iterations and times are not presented.
The stopping criterion was defined by relative residual as $\|r\|_{2}/\|g\|_{2} < 10^{-8}$.

Table~\ref{tab:iter_cavityHh4} presents results for variable size of the problem and variable constraints with fixed ratio $H/h = 4$.
Columns contain results for corners only (`c'), with additional constraints on averages over edges (`c+e'), 
with additional constraints on averages over faces (`c+f'), and with both (`c+e+f').
This experiment shows, that using averages on edges and faces, number of iterations for larger problem does not grow with problem size. 
This is not achieved with corner constraints only, and even with averages on edges or faces, number of iterations slightly grows. 
For reference, number of iterations without preconditioning (`no prec.') is also reported.

Table~\ref{tab:iter_cavity222} presents results for fixed number of subdomains (eight) with variable $H/h$ ratio.
This experiment shows, that number of iterations mildly grows for all types of constraints with growing $H/h$, 
in agreement with available theory for SPD problems.

Numbers of iterations obtained for the same divisions, 
but using Matlab implementation of ILU preconditioner with variable threshold for dropping entries 
(ILUT \cite{Saad-1994-ILU}) are presented in
Tables~\ref{tab:iter_cavityHh4_ilut} and \ref{tab:iter_cavity222_ilut}.
These tables present results with respect to variable threshold for dropping entries in factors ranging from $10^{-3}$ to $10^{-5}$.
Where (`--') occurs, PCG method fails to converge.
We can conclude, that the ILUT preconditioner with threshold $10^{-5}$ is very efficient and seems to reach independence of size of this problem.
These tables also suggest, that with an improving preconditioner, PCG method tends to converge also for these indefinite problems.
Numbers of iterations for threshold $10^{-4}$ are comparable with BDDC with sufficient constraints.

\begin{table}[htbp]
\centering
\begin{tabular}{|c|cccc|c|}
\hline
                  & \multicolumn{4}{|c|}{BDDC}                   & no prec.     \\
subdomains (unknowns)   &  c        &   c+e    &    c+f    &   c+e+f   &              \\ \hline 
8    (15,468)     &  31/24.5  &  28/24.5 &  26/22.5  &  23/19.5  &  223/608.5   \\
27   (49,072)     &  50/74.5  &  35/31.5 &  31/33.5  &  26/23    &  436/731.5   \\
64  (112,724)     & 75/126.5  &  42/46.5 &  34/27.5  &  27/44.5  &  627/1,074   \\
125 (216,024)     & 115/182.5 &  47/41.5 &  35/28.5  &  27/42.5  &  782/1,168.5 \\ \hline
\end{tabular}
\caption{Number of iterations by GMRES/BiCGStab with BDDC preconditioner for variable type of constraints, and without preconditioning, 
variable number of subdomains, $H/h = 4$.}
\label{tab:iter_cavityHh4}
\end{table}

\begin{table}[htbp]
\centering
\begin{tabular}{|c|cccc|c|}
\hline
                  & \multicolumn{4}{|c|}{BDDC}                 &  no prec.  \\
$H/h$ (unknowns)  &  c        &   c+e    &    c+f    &   c+e+f &            \\ \hline 
2   (2,312)       & 26/20     & 22/19.5  & 22/17.5   & 19/15.5 &  92/372.5  \\
4   (15,468)      & 31/24.5   & 28/24.5  & 26/22.5   & 23/19.5 &  223/608.5 \\
8   (112,724)     & 38/44     & 33/26.5  & 29/63.5   & 26/22.5 &  418/690.5 \\
12  (368,572)     & 42/35.5   & 37/31.5  & 32/173.5  & 29/98.5 &  530/725.5 \\ \hline
\end{tabular}
\caption{Number of iterations by GMRES/BiCGStab with BDDC preconditioner for variable type of constraints, and without preconditioning, 
variable $H/h$, 2$\times$2$\times$2 = 8 subdomains.}
\label{tab:iter_cavity222}
\end{table}

\begin{table}[htbp]
\centering
\begin{tabular}{|c|ccc|}
\hline
                  & \multicolumn{3}{|c|}{ILUT threshold}  \\
subdomains (unknowns)   & $10^{-3}$  & $10^{-4}$  & $10^{-5}$   \\ \hline 
8    (15,468)     & 15/11/--   & 5/3.5/--   & 3/2/3       \\
27   (49,072)     & 32/40.5/-- & 9/5.5/--   & 5/2.5/5     \\
64  (112,724)     & 49/97.5/-- & 15/11.5/-- & 6/3.5/20    \\
125 (216,024)     & 57/138/--  & 23/15.5/-- & 6/4/--      \\ \hline
\end{tabular}
\caption{Number of iterations by GMRES/BiCGStab/PCG with ILUT preconditioner for variable threshold, 
variable number of subdomains, $H/h = 4$.}
\label{tab:iter_cavityHh4_ilut}
\end{table}

\begin{table}[htbp]
\centering
\begin{tabular}{|c|ccc|}
\hline
                  & \multicolumn{3}{|c|}{ILUT threshold}  \\
$H/h$ (unknowns)  & $10^{-3}$  & $10^{-4}$  & $10^{-5}$   \\ \hline 
2   (2,312)       & 6/4/11     & 3/1.5/4    & 3/1.5/3     \\
4   (15,468)      & 15/11/--   & 5/3.5/--   & 3/2/3       \\
8   (112,724)     & 43/66.5/-- & 15/10.5/-- & 5/3.5/10    \\
12  (368,572)     & 44/79.5/-- & 29/26/--   & out of memory \\ \hline
\end{tabular}
\caption{Number of iterations by GMRES/BiCGStab/PCG with ILUT preconditioner for variable threshold, 
variable $H/h$, 2$\times$2$\times$2 = 8 subdomains.}
\label{tab:iter_cavity222_ilut}
\end{table}

\subsection{Problem (3)} 

\begin{figure}[htbp]
\begin{center}
   \includegraphics[width=100mm]{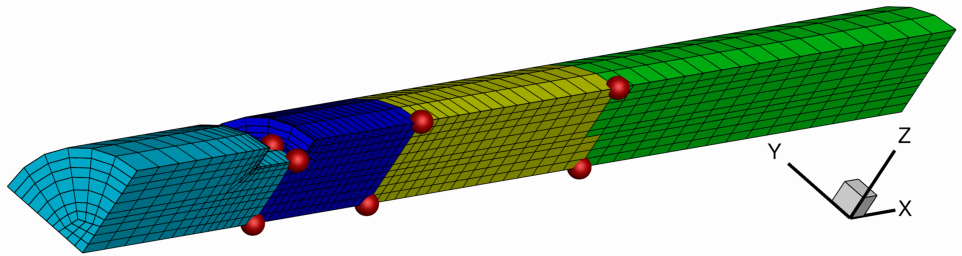}
\end{center}
\caption{\label{fig:tube_mesh} Mesh of the tube of Problem (3) divided into 4 subdomains.}
\end{figure}

\begin{figure}[htbp]
\begin{center}
   \includegraphics[width=65mm]{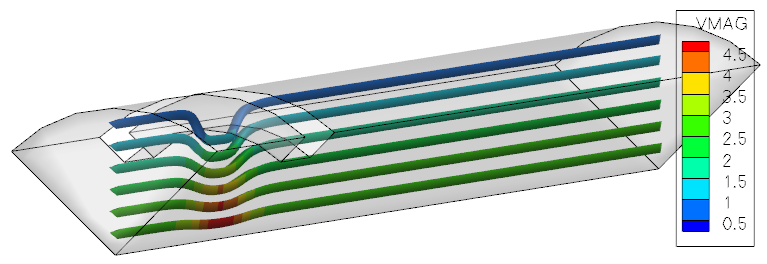}
   \includegraphics[width=65mm]{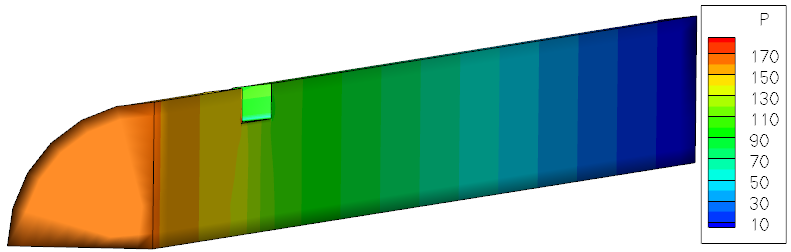} 
\end{center}
\caption{\label{fig:tube_sol} Solution of Problem (3): streamlines coloured by velocity magnitude (left) and pressure contours (right).}
\end{figure}

The last problem is inspired by flow in artificial arteries.
The geometry is simplified to a tube with a~sudden reduction of diameter.
Due to the symmetry of the tube, only one quarter is considered in the computation.
Constant kinematic viscosity $\nu = 0.01$ is considered.
Parabolic velocity profile with unit mean value is prescribed at the inlet,
and `do-nothing' boundary condition (\ref{bcn}) at the outlet.
The diameter of the tube at the inlet is 0.025 and at the narrow part 0.019.
The mesh consists of 3,393 $Q_{2S}Q_{1}$ finite elements with 54,248 unknowns.
It was divided into 4 subdomains by METIS (Figure~\ref{fig:tube_mesh}).
Solution of this rather small problem with only corner constraints requires 33 PCG iterations and takes 30 seconds,
which is comparable to 133 seconds for solution by serial frontal solver, but now obtained in parallel.
Application of averages on faces does not reduce the number of iterations while the solution takes 40 seconds due to the overhead of transforming the matrix.
Streamlines and pressure contours are plotted in Figure~\ref{fig:tube_sol}.


\section{Conclusion}
\label{sec:conclusion}

In our contribution,
we present a straightforward parallel implementation of the BDDC preconditioner based on its global formulation and built on top of the parallel direct solver MUMPS.
After a verification of the solver on a number of problems from linear elasticity analysis,
we explore the application of BDDC to problems with indefinite matrices,
namely the Stokes problem.
Although the available theory either does not cover this case, or treats it differently \cite{Li-2006-BAI,Tu-2005-BAM},
the presented experiments suggest promising ways for this effort.
Results for two versions of the benchmark problem of the lid driven cavity and for a real-life problem are presented. 
These results show that the BDDC preconditioner is applicable to the Stokes flow and may speed up the solution considerably.

Without claiming that this is the general case,
we have performed several experiments, 
for which the current parallel implementation based on the PCG method is successfully used even though the system matrix is indefinite.
The reason why a breakdown was not observed lies probably in the indefiniteness of the BDDC preconditioner.
Although solution times present a~large advancement compared to the method previously used for these problems by our group,
the experiments reveal that optimal scalability is not achieved for the PCG method neither with respect to number of iterations, 
nor with respect to computational times.

On the other hand, our Matlab experiments combining the BDDC preconditioner with GMRES and BiCGStab methods suggest, 
that for suitably chosen constraints, optimal behaviour can be achieved with respect to growing number of subdomains.

\section*{Acknowledgements}
This research has been supported 
by the Czech Science Foundation under Grant GA CR 106/08/0403, 
by National Science Foundation under Grant DMS-0713876,
by the Academy of Sciences of the CR under Grant IAA100760702, 
and by
projects MSM 6840770001 
and MSM 6840770010. 
It has been also supported
by Institutional Research Plan AV0Z10190503.
A part of this work was done during Jakub~\v{S}\'{\i}stek's visit at the University of Colorado Denver.
Jakub~\v{S}\'{\i}stek also thanks to the European Science Foundation OPTPDE (Optimization with PDE Constraints) Network
for the financial support of attending the 10th ICFD Conference on Numerical Methods for Fluid Dynamics.


\end{document}